# The Alternative Form of Fermat's Equation

## Anatoly A. Grinberg

## anatoly_gr@yahoo.com

## ABSTRACT

An alternative form of Fermat's equation[1] is proposed. It represents a portion of the identity that includes three terms of Fermat's original equation. This alternative form permits an elementary and compact proof of the first case of Fermat's Theorem (FT) for a number of specific exponents. Proofs are given for exponents n equal to 3, 5, 7, 11, 13 and 17. All these cases have already been proven using the original Fermat's equation, not to mention the fact that a complete proof of FT was given by A. Wiles [2]. In view of this, the results presented here carry a purely methodological interest. They illustrate the effectiveness and simplicity of the method, compared with the well-known classical approach. An alternative form of the equation permits use of the criterion of the incompatibility of its terms, avoiding the labor-intensive and sophisticated calculations associated with traditional approach.

## I. Introduction

Fermat's last theorem states that there are no three integers A, B and C, which satisfy the equation

$$A^n + B^n - C^n = 0 \qquad\qquad (I.1)$$

when n> 2.

The theorem, formulated by Pierre de Fermat in 1637 [1], attracted the attention of the greatest mathematicians for more than 350 years, in an unsuccessful attempt to find a complete solution.

Efforts to solve FLT led to the development of new directions in mathematics and, together with this, to a significant increase in complexity of the methods used to analyze the problem. This trend culminated in the brilliant achievement of Andrew Wiles, who proved FT in 1994 [2].

All early attempts to solve the problem of Fermat for its particular cases[3,4], were limited to the supposition of the existence of integer solutions of (I.1), with subsequent development of this supposition until it came to contradict some initial condition. Equation (I.1) did not allow for a different approach.

This paper considers an alternative form of the equation, one that is fully equivalent to the original (I.1). Analytically, it is much more complicated than equation (I.1). However, it allows detection of incompatibility of terms contained within it. Inconsistency of terms reflects the absence of solutions in integers for A, B, C. In some cases, to proving incompatibility is much simpler and more productive than a direct attempt to prove the alleged inconsistencies in the hypothetical solution. Moreover, the proof of incompatibility can be made using numerical means, as illustrated below for specific values of n.

The method used in this paper—that of proving incompatibility—is very limited. It is possible that its further development or the formulation of another criterion of incompatibility, will expand the borders of elementary proofs of FT. In this paper we present an elementary proof of the first case FT for prime exponent n equal to 3, 5, 7, 11, 13 and 17.

Let us to describe the well-known characteristics of the variables A, B, C [3, 4]. Without loss of generality, it is assumed that the greatest common divisor of integers A, B and C has been removed. Due to Eq.(I.1), this procedure makes all three numbers relatively prime. This, in turn, leads to that fact that two of the numbers must be odd. Here, as is usually accepted, we will assume that A is even.

## II. Reformulation of Fermat's equation

Consider two binomial identities

$$A^n + B^n - C^n = (A + B - C)^n - \sum_{\nu=1}^{n-1} \binom{n}{\nu} [A^\nu B^{n-\nu} + (A+B)^\nu (-C)^{n-\nu}] \quad \textbf{(II.1)}$$

$$A^n + B^n - C^n = (A + B - C)^n - 2C^n -$$
$$\sum_{\nu=1}^{n-1} \binom{n}{\nu} [A^\nu B^{n-\nu} + (A+B)^\nu (-C)^{n-\nu}] \quad \textbf{(II.2)}$$

The first identity occurs when binomial exponent n is odd; the second identity—for even values of n. Symbol $\binom{n}{\nu}$ denotes the binomial coefficients.

Both identities are valid for an arbitrary number of A, B and C. On the left, in (II.1) and (II.2), there are three terms of equation (I.1). Equating their sum to zero, we thus equate to zero the right-hand sides of the identities, in order to give another equivalent form of Fermat's equation. Naturally, the values of A, B and C, in which the right and left sides of the identities (II.1) and (II.2) vanish, must also be identical. Or, closer to the subject, if the left side of identities (II.1) and (II.2) does not vanish for any integer values of the variables A, B and C, then the same should hold for the right-hand sides of these identities.

In this new formulation, equation (I.1) splits into two equations

$$(A + B - C)^n = \sum_{\nu=1}^{n-1} \binom{n}{\nu} A^{\nu} B^{n-\nu} + \sum_{\nu=1}^{n-1} \binom{n}{\nu} (A + B)^{\nu} (-C)^{n-\nu} \qquad \textbf{(II.3)}$$

$$(A + B - C)^n = 2C^n + \sum_{\nu=1}^{n-1} \binom{n}{\nu} A^{\nu} B^{n-\nu} + \sum_{\nu=1}^{n-1} \binom{n}{\nu} (A + B)^{\nu} (-C)^{n-\nu} \qquad \textbf{(II.4)}$$

for odd and even exponent n, respectively.

When using these equations as an alternative to the Fermat equation (I.1), we should not use the latter, in order to avoid coming to the identity A = A.

Two sums on the right-hand sides of these equations are truncated binomials of degree n in two pairs of numbers: (A, B), and (A + B, -C). Generally, a truncated binomial U (a, b) from the sum of the two numbers (a, b) we defined as

$$\textbf{U}(\boldsymbol{a}, \boldsymbol{b}) = (\boldsymbol{a} + \boldsymbol{b})^n - \boldsymbol{a}^n - \boldsymbol{b}^n =$$

$$\sum_{\nu=1}^{n-1} \binom{n}{\nu} \boldsymbol{a}^{\nu} \boldsymbol{b}^{n-\nu} = -\sum_{\nu=1}^{n-1} \binom{n}{\nu} \boldsymbol{q}^{\nu}(-\boldsymbol{b})^{n-\nu} \qquad \textbf{(II.5)}$$

where $\quad \textbf{q} = \textbf{a} + \textbf{b.}$

In other words, from $U(a, b)$ are removed the side terms of a full binomial series. We confine ourselves to the cases where n is a prime number. In this case, all the binomial coefficients in (II.5) are proportional to n. Furthermore, regardless of the parity of n and integers $a$ and $b$, the truncated series $U(a, b)$ is always an even number.

Apart from the obvious divisibility of $U(a, b)$ by $a$, b and $q = (a + b)$, this series is also divisible by $(a^2 + ab + b^2)$ for n = 6k +1 and n = 6k +5. Moreover, at n=6k+1, the series $U(a, b)$ is divisible by the square of $(a^2 + ab + b^2)$ [8].

From Eq.(II.5) it follows [5], that, if only one of the numbers $a, b$ or $q = (a + b)$ is a multiple of n, then $U(a, b)$ is divisible only by $n^2$.

By virtue of the described properties of the truncated binomial, the left side of equation (II.3) must be divisible by 2n. Since n, as well as 2, are prime integers, (A + B - C) is divisible by 2n. It can, therefore, be represented as

$$A + B - C = 2\beta n \qquad \textbf{(II.6)}$$

where β is the integer coefficient.

Substituting (II.6) in equation (II.3), we obtain

$$(2\,\beta\,n)^n = U(A,B) + U(A+B,-C) \qquad \textbf{(II.7)}$$

The left side of equation (II.7) imposes strict restrictions on the right side: the latter should be divisible by left-standing factors, in particular, by $n^n$. If one can show that division by at least one of these factors is impossible, it would show the contradictory nature of the equation. Inconsistency, of course, arises from the fact that, in agreement with Fermat's statement, the right side of identity (II.1) cannot become zero under any integer values of the variables A, B and C. However, this condition is necessary but not sufficient for Eq.(II.7) to exist in integers[7].

It is easy to see that division by $n^n$ might be the case if all the variables A, B and C are divisible by n. But this would contradict the mutual primacy of these integers. We restrict ourselves to the first case of TF, in which none of the variables A, B, C is divisible by n.

By condition (II.6), the sum of the arguments of a truncated series of U (A + B, -C) is a multiple of n. As mentioned above, in this case, U(A+ B, -C) is divisible only by $n^2$. Therefore, if U (A, B) is divided by $n^k$, for k $\neq$ 2, then equation (II.7) does not exist in integers. A more detailed analysis of compatibility of equations of type (II.7) is given in [7]. Below we give elementary proofs of the first case TF for n equal to 3, 5, 7, 11, 13 and 17. All these cases were proved by classical methods using the original form (I.1) of Fermat's equation.

## 1. Proof of TF for n = 3.

For n = 3, the truncated series in Eq. (II.7) are determined by expressions

$$U(A,B) = n\,A\,B\,(A+B) \qquad \textbf{(II.8a)}$$

$$U(A+B,-C) = -n(A+B)\,C\,(A+B-C) = -2\beta n^2\,(A+B)C \qquad \textbf{(II.8b)}$$

so that Eq. (II.7) takes the form

$$(2\,\beta\,n)^n = n\,(A+B)\,(A\,B - 2\beta n\,C\,) , \qquad \textit{(II.9)}$$

Since C is not divisible by n, then, according to (II.6), (A + B) also is not divisible by n. Hence, the right side of equation (II.9) is a multiple of only n. This proves FT for n = 3.

## 2. Proof of TF for n = 5.

For n = 5 we find

$$\mathbf{U}(A, B) = \boldsymbol{n\,A\,B\,(A + B)\,(A^2 + A\,B + B^2)} \qquad \textbf{(II.10a)}$$

$$\mathbf{U(A + B, -C)} = -\,2\boldsymbol{\beta n^2}(A + B)\,\mathbf{C}\,[(2\boldsymbol{\beta n})^2 +\ 2\boldsymbol{\beta n}\,\mathbf{C} + \mathbf{C^2}] \qquad \textbf{(II.10b)}$$

It follows from the general analysis that U(A + B,-C) is divisible by $n^2$. The divisibility by n of U(A, B)/n is determined by the divisibility of factor $(A^2 + AB + B^2)$. To find out whether this factor is divisible by n, it is convenient to use modular arithmetic. We then have the following comparison

$$(\mathbf{A^2 + AB + B^2}) \equiv \ \Delta_\mathbf{A}{}^2 + \Delta_\mathbf{A}\,\Delta_\mathbf{B} + +\Delta_\mathbf{B}{}^2 \quad (\boldsymbol{mod}\ \mathbf{5}) \qquad \textbf{(II.11)}$$

where $\Delta_A$ and $\Delta_B$ represent the remainders of numbers A and B, respectively, modulo n=5. Since A and B are not divisible by n, the values of $\Delta_A$ and $\Delta_B$ are limited to a set of numbers {1 , 2 , . . . (n-1)}. Thus, there $are$ $(n-1)^2$ pairs of combinations of numbers $\Delta_A$, $\Delta_B$ for which the value of

$$\mathbf{F}(\Delta_\mathbf{A}, \Delta_\mathbf{B}) = \ \Delta_\mathbf{A}{}^2 + \Delta_\mathbf{A}\,\Delta_\mathbf{B} + \Delta_\mathbf{B}{}^2 \qquad \textbf{(II.12)}$$

needs to be tested for its multiplicity of n. Function $F(\Delta_A, \Delta_B)$ is symmetric with respect to its arguments. Therefore, it is sufficient to test only those combinations for which $\Delta_A \geq \Delta_B$. The number of such combinations is equal to $\frac{1}{2}$ [n(n − 1) ]. Because the sum (A + B) is also not divisible by n, the combinations for which $(\Delta_A +\ \Delta_B) = n$, can be excluded. The number of combinations satisfying conditions $(\Delta_A +\ \Delta_B) = n$ and $\Delta_A \geq \Delta_B$, equals (n - 1) / 2. Thus, the total number of combinations $\Delta_A$, $\Delta_B$ that must be tested is equal to $N = (n-1)^2/2$. For n = 5 , none of the eight (N=8) such combinations does not lead to $F(\Delta_A, \Delta_B)$ multiple of n. This proves the TF for n = 5

We note that the above expression of the number of combinations N that must be tested, is valid for any algebraic function F(A,B) that is symmetrical with respect to its arguments. The truncated binomial U(A,B) is an example of such a function.

## 3. Proof of TF for n = 7

At n = 7, equation (II.7) takes the form

$$(2\boldsymbol{\beta n})^n = n(A + B)[AB(A^2 + AB + B^2\,)^2 + 2\boldsymbol{\beta Cn}((A + B)^2 + 2\boldsymbol{\beta Cn})^2] \quad \textbf{(II.13)}$$

from which it follows that the proof of the incompatibility of its terms is not required for any calculations. Indeed, if the trinomial $(A^2 + AB + B^2)$ is not divisible by n=7 then the right side of Eq. (II.13) is divisible only by first degree n . On other hand, if this trinomial is divisible by n, then divisibility of the right side of Eq. (II.13) is determined by the second term in square brackets. Therefore, in this case, the right side of Eq. (II.13) is only divisible by the square of n. So, in both cases, it is not divisible by $n^7$, which proves the FT for n = 7.

## 4. Proof of FLT for n = 11

The term U(A, B) in equation (II.7), for n=11, is determined by the expression

$$U(A, B) = n\,A\,B\,(A + B)(A^2 + AB + B^2)\,(A^6 + 3A^5B +$$

$$7A^4B^2 + 9A^3B^3 + 7A^2B^4 + 3AB^5 + B^6\,) \qquad \textbf{(II.14)}$$

The divisibility of ratio U(A, B)/n by n was tested numerically by a method similar to that used in the case of n = 5. Calculations have shown that this ratio is not divisible by n = 11, and hence the terms of equation (II.7) are incompatible.

## 5. Proof of FLT for n = 13 and 17

When n = 13, U (a, b) in equation (II.7) has the form

$$U(A, B) = n\,A\,B\,(A + B)(A^2 + AB + B^2)^2$$

$$\left(A^6 + 3A^5B + 8A^4B^2 + 11A^3B^3 + 8A^2B^4 + 3AB^5 + B^6\right) \quad \textbf{(II.15)}$$

Numerical analysis similar to that used in the description of the cases n = 5 and 11, shows that the last factor in (II.15) is not divisible by n = 13, while the trinomial ($A^2 + AB + B^2$) is divisible by n for some values of A and B. However, as in the case of n = 7, this only leads to the fact that the divisibility of the right side of equation (II.7) is determined by the divisibility of truncated binomial U(A + B, -C), which is only divisible by $n^2$.

Using similar approach it could be shown that at n=17 the truncated binomial $U(A, B)$ is not divisible by $n^2$.

We note here that, for the second case of FT, when only one of the numbers A, B, C is a multiple of n, the method of proof proposed here is not applicable [7].

## III. Conclusion

The afore-presented method of analysis of Fermat's equation is not unique. Here we used a method based on the analysis of the compatibility of the various terms of the equation. Despite its simplicity, or, perhaps, because of it, it is not obvious that a more advanced elementary proof of FT can be realized by following only this approach. I want to reiterate that this is a simple method for proving FT, and that I remain mindful of the fact that the full, complete proof of Fermat's Last Theorem was implemented by Andrew Wiles in 1994. [2]